\documentclass[12pt]{amsart}
\title[Co-Hopfian Sierpi\'nski carpet]{A Sierpi\'nski carpet with\\ the co-Hopfian property
}
\author{Sergei Merenkov}
\address{Department of Mathematics\\
University of Illinois\\ 1409 W. Green Street\\ Urbana, IL
61801\\USA} \email{merenkov@illinois.edu}
\thanks{
Supported by NSF grant
DMS-0653439.}

\usepackage[active]{srcltx}

\usepackage{amsmath, amssymb, amsthm,latexsym}
\usepackage{graphicx}

\newcommand\N{{\mathbb N}}

\newcommand\Z{{\mathbb Z}}
\newcommand\R{{\mathbb R}}

\newcommand\mH{{\mathcal H}}

\newcommand\dee{\partial}
\newcommand\diam{\operatorname{diam}}
\newcommand\id{\operatorname{id}}

\renewcommand\:{\colon}

\renewcommand\Im{\mathop{\mathrm{Im}}}

\newcommand\no{\noindent}

\newtheorem{theorem}{Theorem}[section]

\newtheorem{proposition}[theorem]{Proposition}
\newtheorem{corollary}[theorem]{Corollary}

\newtheorem{lemma}[theorem]{Lemma}

\theoremstyle{definition}

\begin{document}

\abstract{Motivated by questions in geometric group theory we define a quasisymmetric co-Hopfian property for metric spaces and provide an example of a metric Sierpi\'nski carpet with this property. As an application we obtain a quasi-isometrically co-Hopfian Gromov hyperbolic space with a Sierpi\'nski carpet boundary at infinity. In addition, we give a complete description of the quasisymmetry group of the constructed Sierpi\'nski carpet. This group is uncountable and coincides with the group of bi-Lipschitz transformations.
}
\endabstract

\maketitle

\section{Introduction}\label{S:Intro}

\no
The co-Hopfian property studied in geometric group theory is defined as follows. A group $G$ is said to be \emph{co-Hopfian} if every monomorphism of $G$ into itself is an isomorphism. Every finite group is obviously co-Hopfian and many  
examples of infinite groups possessing this property can be found in~\cite{FH07}, \cite{pH00}, \cite{KW01}, \cite{zS97}.
A related co-Hopfian property for unbounded metric spaces
is defined as follows. A  map $\phi$ of a metric space $(X, d_X)$ to a metric space $(Y, d_Y)$ is a \emph{quasi-isometric embedding} if there exist constants $\lambda\geq1$ and $C\geq 0$ such that 
$$
\frac1{\lambda}d_X(p,q)-C\leq d_Y(\phi(p),\phi(q))\leq \lambda \cdot d_X(p,q)+C
$$
for all $p, q\in X$. If we want to emphasize the parameters $\lambda$ and $C$, we say that $\phi$ is a $(\lambda, C)$-\emph{quasi-isometric embedding}. A quasi-isometric embedding of $(X, d_X)$ to $(Y, d_Y)$ is called a \emph{quasi-isometry} if in addition there exists a constant $D\geq0$ such that every point in $Y$ is within distance $D$ from $\phi(X)$. The two spaces $(X, d_X)$ and $(Y, d_Y)$ are then called \emph{quasi-isometric}, and
this is an equivalence relation.  
We say that a metric space $(X,d)$ is \emph{quasi-isometrically co-Hopfian} if every quasi-isometric embedding of $X$ into itself is a quasi-isometry. 

It is known that uniformly contractible, bounded geometry manifolds, e.g., Euclidean spaces, are quasi-isometrically co-Hopfian, as are coarse $PD(n)$ spaces, see~\cite{KK05}.
Many quasi-isometrically co-Hopfian metric spaces
can be found among Gromov hyperbolic spaces, see~\cite{BS00}, \cite{BS07} for background and terminology on general Gromov hyperbolic spaces and~\cite{GH90}, \cite{mG87} for background on Gromov hyperbolic groups. 

In this paper we consider only roughly geodesic Gromov hyperbolic spaces. If $(X,d)$ is a metric space, a \emph{roughly quasi-isometric path} in $X$ is a $(\lambda, C)$-quasi-isometric embedding of a segment $[a,b]\subset\R$ or a half-open interval $[a,b)\subset\R$ to $X$. In the latter case we may assume $b=\infty$. 
A $C$-\emph{roughly geodesic segment} $[p,q]$ in $X$ is a $(1,C)$-quasi-isometric embedding $i$ of a segment $[a,b]\subset \R$ into $X$ such that $i(a)=p, i(b)=q$.  A $C$-\emph{roughly geodesic ray} in $X$ is a $(1,C)$-quasi-isometric embedding of $[a,\infty)$ to $X$ for some $a\in\R$. A metric space $(X,d)$ is said to be \emph{roughly geodesic} if there exists $C\geq 0$ such that for every two points $p$ and $q$ in $X$ there is a $C$-{roughly geodesic segment} $[p,q]$. In a roughly geodesic metric space one can speak of \emph{roughly geodesic triangles} $[p,q]\cup[q,q']\cup[q',p]$ formed by roughly geodesic segments $[p,q], [q,q']$, and $[q',p]$. A roughly geodesic metric space is \emph{Gromov hyperbolic} if there exists $\delta\geq0$ such that every roughly geodesic triangle is $\delta$-\emph{thin}. The latter means that every side of such a triangle is contained in the $\delta$-neighborhood of the union of the other two sides. This is a version of the Rips thin triangles definition of Gromov hyperbolicity and it is an exercise to check that for roughly geodesic metric spaces this  notion of hyperbolicity agrees with the one defined via the Gromov product.
Every roughly geodesic Gromov hyperbolic space $X$ has an associated notion of the \emph{boundary at infinity} $\dee_\infty X$. The elements of $\dee_\infty X$ are equivalence classes of roughly geodesic rays, where two such rays are \emph{equivalent} if the Hausdorff distance between them is finite. If $\gamma$ is a roughly geodesic ray in $X$, we refer to the equivalence class containing $\gamma$ as its \emph{endpoint}. If $X$ is a Gromov hyperbolic space, its boundary at infinity $\dee_\infty X$ carries a canonical family of so-called ``visual" metrics, see~\cite[Lemma~6.1]{BS00}, 
and in each of these metrics $\dee_\infty X$ is bounded and complete, see~\cite[Proposition~6.2]{BS00}. 
In what follows it is assumed that  roughly geodesic Gromov hyperbolic spaces are \emph{visual}, i.e., there exists a base point $o$ in $X$ and a constant $C\geq 0$ such that every point $p\in X$ lies on a $C$-roughly geodesic ray $i\: [a,\infty)\to X$ with $i(a)=o$. It is an elementary fact that if $X$ is visual with respect to some base point $o$, then it is visual with respect to any other base point $o'$. It is also trivial that 
for every $p\in X$  and every element $\xi\in\dee_\infty X$ there exists a roughly geodesic ray $i\:[a,\infty)\to X$ in $\xi$ such that $i(a)=p$.

The quasi-isometric co-Hopfian property for Gromov hyperbolic spaces is closely related to the so-called quasisymmetric co-Hopfian property of their boundaries at infinity. Recall that a homeomorphism $f$ between metric spaces $(X,d_X)$ and $(Y,d_Y)$ is called \emph{quasisymmetric} 
if there exists a homeomorphism $\eta\: [0,\infty)\to[0,\infty)$ such that 
$$
\frac{d_Y(f(p),f(q))}{d_Y(f(p),f(q'))}\leq\eta\bigg(\frac{d_X(p,q)}{d_X(p,q')}\bigg)
$$ 
for all triples of distinct points $p,q$, and $q'$ in $X$. 
A \emph{quasisymmetric embedding} $f$ of $(X,d_X)$ into $(Y,d_Y)$ is a one-to-one continuous map of $X$ into $Y$ that is quasisymmetric between $(X,d_X)$ and $(f(X),d_Y)$.  We say that a metric space $(X,d_X)$ is \emph{quasisymmetrically co-Hopfian} if every quasisymmetric embedding of $X$ into itself is onto.
 
If  $(X,d)$ is a visual roughly geodesic Gromov hyperbolic space such that $\dee_\infty X$ is quasisymmetrically co-Hopfian, then $X$ is quasi-isometrically co-Hopfian. 
Indeed, every quasi-isometric embedding $\phi$ of a roughly geodesic Gromov hyperbolic space $X$ into itself induces a quasisymmetric embedding $\dee\phi$ of $\dee_\infty X$ into itself, see~\cite[Theorem~6.5]{BS00}. 
If $\dee X_\infty$ is quasisymmetrically co-Hopfian, the map $\dee\phi$ is onto.
Now let $p$ be a point in $X$. Since $X$ is visual, $p$ lies on a roughly geodesic ray $\tilde \gamma$ emanating from some base point $o\in\phi(X)$. Let $\tilde \xi\in\dee_\infty X$ denote the endpoint of $\tilde\gamma$ and $\xi=(\dee\phi)^{-1}(\tilde\xi)$. If $\gamma\: [a,\infty)\to X$ is a roughly geodesic ray in $\xi$ such that $\gamma(a)\in\phi^{-1}(o)$, then $\phi(\gamma)\: [a,\infty)\to X$ is a roughly quasi-isometric path with $\phi(\gamma(a))=o$. The path $\phi(\gamma)$ is not necessarily a roughly geodesic ray, but $\phi(\gamma)$ still has a well-defined notion of an endpoint in $\dee_\infty X$ and in our case the endpoint is $\tilde \xi$, see~\cite[Proposition~6.3]{BS00}. Applying the stability of roughly quasi-isometric paths~\cite[Proposition~5.4]{BS00} to $\tilde\gamma$ and $\phi(\gamma)$, we conclude that $p$ is within bounded distance from $\phi(X)$, i.e., $\phi$ is a quasi-isometry. 

A large source of examples of Gromov hyperbolic spaces is the theory of Gromov hyperbolic groups. If $G$ is a finitely generated group and $S$ is a finite symmetric (i.e., $S$ contains the inverse of each of its elements) set of generators, 
then one can consider the \emph{Cayley graph} $\Gamma(G,S)$. The vertices of $\Gamma(G,S)$ are the elements of the group $G$ and two vertices $v_1$ and $v_2$ are connected if and only if $v_1^{-1}v_2\in S$. The Cayley graph can be made into a geodesic metric space in a natural way by declaring each edge with the two vertices as its endpoints to be isometric to the unit segment $[0,1]$. 
A finitely generated group $G$ is called \emph{Gromov hyperbolic} if for some finite symmetric generating set $S$, the metric space $\Gamma(G,S)$ is Gromov hyperbolic. 
It is a fact that if $S_1$ and $S_2$ are two finite symmetric generating sets, then $\Gamma(G,S_1)$ is Gromov hyperbolic if and only if $\Gamma(G, S_2)$ is Gromov hyperbolic. 

\begin{figure}
[htbp]
\begin{center}
\includegraphics[height=40mm]{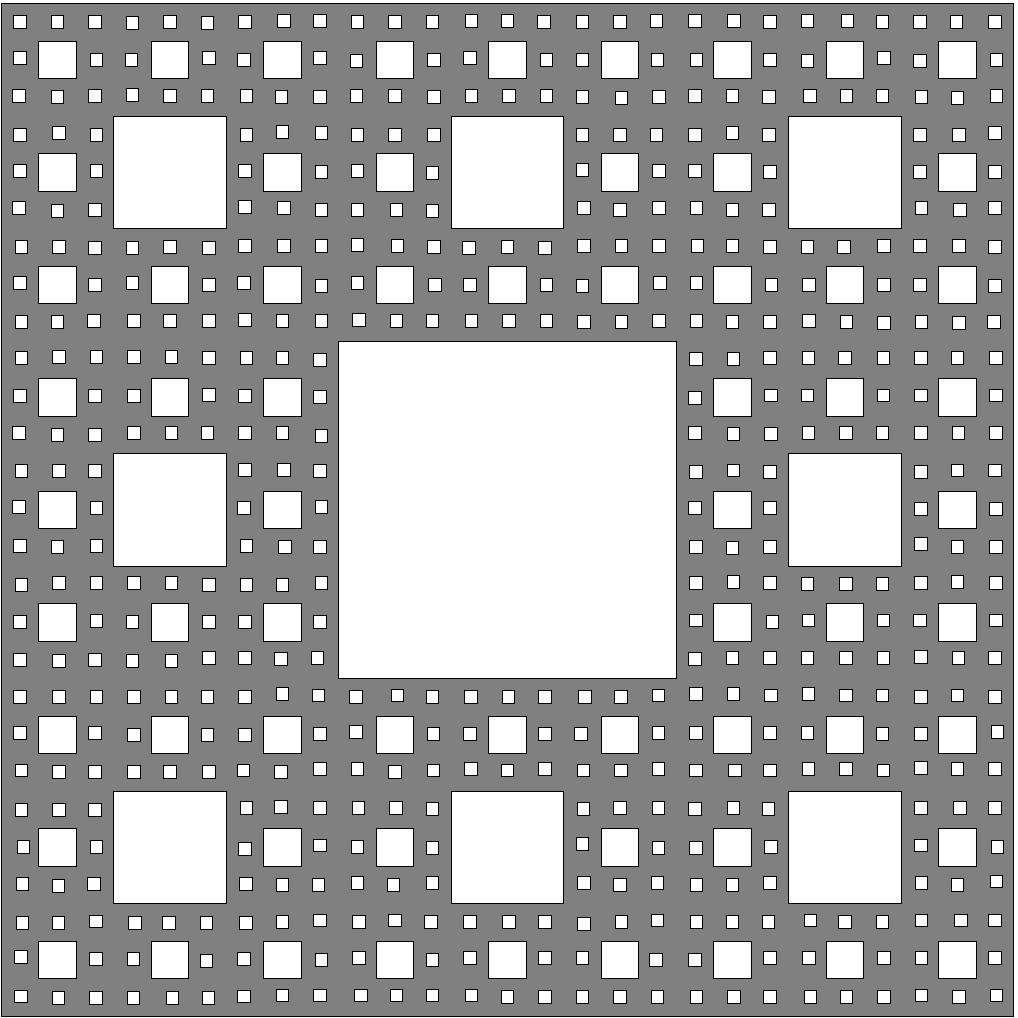}
\caption{
The standard Sierpi\'nski carpet $S_3$.
}
\label{F:St}
\end{center}
\end{figure}

If $G$ is the fundamental group of a closed hyperbolic manifold, it is Gromov hyperbolic and its boundary at infinity $\dee_\infty G$ is a topological sphere. Topological spheres have a stronger co-Hopfian property, namely every continuous embedding into itself is onto.
We call topological spaces with such a property  \emph{topologically co-Hopfian}. Clearly, every topologically co-Hopfian metric space is quasisymetrically co-Hopfian, and therefore the fundamental group of every closed hyperbolic manifold is quasi-isometrically co-Hopfian.
There are easy examples of unbounded or non-complete spaces that are quasisymmetrically co-Hopfian but not topologically co-Hopfian, e.g., Euclidean spaces or standard spheres with finitely many punctures. 
We are mostly interested however in the co-Hopfian property of bounded complete metric spaces since boundaries at infinity are such spaces. 
The primary purpose of this paper is to provide an example of a Gromov hyperbolic space that is quasi-isometrically co-Hopfian for non-topological reasons. 
\begin{theorem}\label{T:QiCoHopf}
There exists a quasi-isometrically co-Hopfian visual roughly geodesic Gromov hyperbolic space $X$ whose boundary at infinity $\dee_\infty X$ is a Sierpi\'nski carpet.
\end{theorem}
A \emph{Sierpi\'nski carpet} is a compact topological space homeomorphic to the standard Sierpi\'nski carpet $S_3$, see Figure~\ref{F:St}. 

Given any bounded complete metric space $(Z, d_Z)$, there is a  visual roughly geodesic Gromov hyperbolic metric space, called the ``cone'' of $Z$ and denoted Con$(Z)$, whose boundary at infinity is $Z$ and $d_Z$ is bi-Lipschitz to a visual metric~\cite[Theorems~7.2, 8.1]{BS00}.
Therefore, Theorem~\ref{T:QiCoHopf} follows from the following result.
\begin{theorem}\label{T:QsCoHopf}
There exists a metric Sierpi\'nski carpet that is quasisymmetrically co-Hopfian. 
\end{theorem}
The Sierpi\'nski carpet in Theorem~\ref{T:QsCoHopf}
is a double of a self-similar Sierpi\'nski carpet 
and it has many nice geometric and analytic properties: 
its \emph{peripheral circles}, i.e., embedded non-separating simple closed curves, are uniformly relatively separated uniform quasicircles, it is linearly locally connected, Ahlfors 2-regular in the Hausdorff 2-measure, and therefore doubling. However, it is not Loewner and does not satisfy a (1,2)-Poincar\'e inequality, see~\cite{BK02}, \cite{jH01} for the definitions.

Not every bounded complete metric space can arise as the boundary at infinity of a Gromov hyperbolic group. For example, if $G$ is a Gromov hyperbolic group such that its boundary at infinity $\dee_\infty G$ has a manifold point, then $\dee_\infty G$ must be a topological sphere~\cite[Theorem 4.4]{KB02}.
Such a group is then quasi-isometrically co-Hopfian. In general it is hard to establish the quasi-isometric co-Hopfian property for groups. 
For example, it is unknown whether the fundamental groups of compact hyperbolic 3-manifolds with non-empty totally geodesic boundaries are quasi-isometrically co-Hopfian. The boundaries at infinity for these groups are {Sierpi\'nski carpets}.

Sierpi\'nski carpets are in a sense the simplest connected non-manifold boundaries of Gromov hyperbolic groups.  
Indeed, if $G$ is a Gromov hyperbolic group that does not split over a finite or a virtually cyclic group and the boundary at infinity $\dee_\infty G$ has topological dimension one, then $\dee_\infty G$ is homeomorphic to either a circle, or the standard Sierpi\'nski carpet, or the Menger curve~\cite[Theorem 1]{KK00}. Groups with circle boundaries are well understood and are co-Hopfian. The uniformization of Gromov hyperbolic groups whose boundaries at infinity are Sierpi\'nski carpets is addressed by the Kapovich--Kleiner conjecture~\cite{KK00} and it is unknown whether such a group can be co-Hopfian. Groups whose boundaries at infinity are homeomorphic to the Menger curve are generic and the question whether such groups can or cannot have the co-Hopfian property is widely open.

The quasisymmetric co-Hopfian property of many other interesting compact metric spaces is either false or unknown. If a compact manifold has a boundary point it cannot be quasisymmetrically co-Hopfian. This can be seen by pushing the boundary inside the manifold locally near a boundary point, and it can be done quasisymmetrically. 
The standard Sierpi\'nski carpet $S_3$ with the restriction of the Euclidean metric is not quasisymmetrically co-Hopfian since it is metrically self-similar. 
It is an open question whether 
a double of $S_3$ across a peripheral circle is co-Hopfian.
It is unknown if there are metric spaces homeomorphic to the Menger curve, in particular the boundaries of
the Bourdon--Pajot hyperbolic buildings~\cite{BP00}, that are quasisymmetrically co-Hopfian.

\smallskip
\noindent
{\bf Acknowledgment.} The author is grateful to  Ilya Kapovich for suggesting the problem of finding metric spaces that are co-Hopfian for non-topological reasons. He thanks Mario Bonk for explaining various properties of slit carpets used in the present construction. He thanks Kevin Pilgrim for renewing interest in slit carpets. He thanks John Mackay and Jeremy Tyson for many interesting discussions. 
He thanks Maria Sabitova for many stimulating conversations during the revision stage of the paper. 
Last but not the least, the author would like to thank the anonymous referee for numerous comments and suggestions. 

\section{Slit carpets}

\begin{figure}
[htbp]
\begin{center}
\includegraphics[height=40mm]{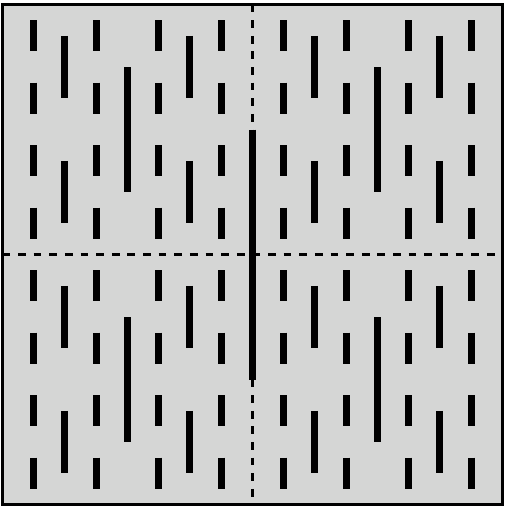}
\caption{
Slit carpet $S_2$.
}
\label{F:Slitc}
\end{center}
\end{figure}

\no
Consider the following space, denoted by $S_2$. We start with the closed unit square $[0,1]\times[0,1]$ in the plane and subdivide it into subsquares
$$
[0,1/2]\times[1/2,1],\ [1/2,1]\times[1/2,1],\ [0,1/2]\times[0,1/2],\ [1/2,1]\times[0,1/2].
$$
We then slit it in the vertical interval connecting the points $(1/2, 1/4)$ and $(1/2, 3/4)$, i.e., double each point of this open interval, and apply such rescaled operations on the four subsquares. This process is continued indefinitely, see Figure~\ref{F:Slitc}. 

Formally, the space $S_2$ can be defined as follows. Let $Q_n,\ n\geq1$, be the finitely connected domain obtained from the open unit square $Q_0=(0,1)\times(0,1)$ by removing the closures of all the slits in the construction of $S_2$ up to the $n$'th generation, i.e., all the slits whose length is at least $1/2^n$. We denote by $\bar Q_n,\ n\geq0$, the completion of $Q_n$ in the path metric $d_{\bar Q_n}$ induced by the Euclidean metric in the plane, and we call the boundary components of $\bar Q_n$ that correspond to the slits of $S_2$ \emph{slits} and the remaining boundary component the \emph{outer square}. 
For every $m,n\in\N\cup\{0\}$ with $m\leq n$ there is a natural 1-Lipschitz projection $\pi_{mn}\: \bar Q_n\to\bar Q_m$ obtained by identifying the points on the slits of $\bar Q_n$ that correspond to the same point of $\bar Q_m$.

As a topological space $S_2$ is the inverse limit of the system $(\bar Q_n, \pi_{mn})$, and as such it is a compact Hausdorff space. For each $n\in\N\cup\{0\}$, the natural projection of $S_2$ onto $\bar Q_n$ will be denoted by $\pi_n$.
The \emph{slits} and the \emph{outer square} of $S_2$ are the topological circles that are the inverse limits of  the slits and the outer squares of $\bar Q_n,\ n\geq 1$, respectively. Clearly, the slits are dense in $S_2$, i.e., for every point $p$ in $S_2$ and every neighborhood $U$ of $p$, there exists a slit of $S_2$ that intersects $U$. 

The diameter of each $\bar Q_n$ is clearly bounded by 3.
If $p=(p_0, p_1,\dots)$, $q=(q_0, q_2,\dots)\in S_2$, we define a distance between them by 
$$
d_{S_2}(p,q)=\lim d_{\bar Q_n}(p_n,q_n).
$$
Since every $\pi_{mn}$ is 1-Lipschitz, $(d_{\bar Q_n}(p_n,q_n))$ is a monotone increasing bounded sequence, and thus $d_{S_2}(p,q)$ exists and defines a metric on $S_2$. 

A \emph{curve} in a metric space $(X,d)$ is a continuous map of one of the intervals $[a,b], [a,b), (a,b]$, or $(a,b)$ into $X$. 
A curve $\gamma$ is said to be \emph{geodesic} if the distance between any two points $p$ and $q$ on $\gamma$ is equal to the length of the part of $\gamma$ between $p$ and $q$. A metric space $(X,d)$ is called \emph{geodesic} if any two points in $X$ can be connected by a geodesic.
Each space $(\bar Q_n, d_{\bar Q_n}),\ n\geq0$, is clearly {geodesic}. The metric $d_{\bar Q_n}$ being  the path metric induced by the Euclidean metric in the plane means that the length of any curve $\gamma$ in $\bar Q_n$ is equal to the Euclidean length of $\pi_{0n}(\gamma)$. 
If $p$ and $q$ are two points in $S_2$ and $n\in\N$, let $\gamma_n$ be an arc-length parametrized geodesic in $\bar Q_n$ connecting $p_n$ and $q_n$. Then 
$$
\tilde\gamma_{n}=\pi_{0n}\circ\gamma_n\circ\frac{d_{\bar Q_n}(p_n,q_n)}{d_{S_2}(p,q)}\: [0,d_{S_2}(p,q)]\to\bar Q_0
$$
is a 1-Lipschitz map for every $n$.
Using the Arzel\`a--Ascoli theorem we conclude that there exists a subsequence of $(\tilde \gamma_{n})$ that converges uniformly to a curve $\tilde\gamma$ in $\bar Q_0$. The map $t\mapsto\tilde\gamma(t)$ is 1-Lipschitz from $[0,d_{S_2}(p,q)]$ to $\bar Q_0$. From the definition of $\tilde\gamma$ we see that $\tilde\gamma$ \emph{lifts} to each $\bar Q_n,\ n\geq1$, i.e., there exists a curve $\gamma_n'$ in $\bar Q_n$ such that $\pi_{0n}(\gamma_n')=\tilde\gamma$, and from the universality property of the inverse limit it follows that $\tilde\gamma$ lifts to a curve $\gamma$ in $S_2$ that connects $p$ and $q$. The length of $\gamma$ is thus at least $d_{S_2}(p,q)$ and it cannot be larger than $d_{S_2}(p,q)$ because $\tilde\gamma$ is 1-Lipschitz. 
This readily implies that $S_2$ is a path-connected geodesic metric space and $d_{S_2}$ is the path metric on $S_2$ induced by the Euclidean metric in the plane. 
It is a simple exercise to check that the topology defined by this metric agrees with the topology of the inverse limit. For each $n\in\N\cup\{0\}$, the natural projection $\pi_n$ of $S_2$ onto $\bar Q_n$ is 1-Lipschitz. To simplify the notation below we denote $\pi_0$ by $\pi$. 
\begin{lemma}\label{L:Sierp}
The space $S_2$ is a Sierpi\'nski carpet whose peripheral circles are the slits along with the outer square. 
\end{lemma}
\noindent
\emph{Proof.}
To show that $S_2$ is a Sierpi\'nski carpet we find a Lipschitz embedding of this space into $\R^2$ and check that the image $S$ is a set that is obtained from the closure of a Jordan domain $D_0$ by removing a countable collection of Jordan domains $D_n,\ n\geq 1$,  so that the following properties are satisfied. The boundaries $\dee D_n,\ n\geq 0$, are pairwise disjoint, they form a null sequence, i.e., ${\rm diam}(\dee D_n)\to0$, and the remaining set $S=\bar D_0\setminus\cup_{n\geq1}D_n$ has no interior. Whyburn's characterization~\cite{gW58} then gives that $S_2$ is a Sierpi\'nski carpet and its peripheral circles are the preimages of $\dee D_n,\ n\geq0$, under the embedding.

A Lipschitz embedding can be obtained inductively as follows. Clearly, there is a $C_1$-Lipschitz embedding $L_1$ of $\bar Q_1$ into $\bar Q_0\subset\R^2$ with $C_1>1$ arbitrarily close to 1, and it agrees with $\pi_{01}$ on the outer square. Geometrically, $L_1$ is obtained by ``opening the slit slightly''. 
Assume that there is a $C_n$-Lipschitz embedding $L_n$ of $\bar Q_n$ into $\bar Q_0$. 
For all $m,n\in\N\cup\{0\}$ with $m\leq n$, there is a natural partition, denoted by $\mathcal P_{mn}$,
of $\bar Q_n$ by the rescaled copies of $\bar Q_m$.
Every rescaled copy of $\bar Q_1$ in $\mathcal P_{1(n+1)}$ is mapped by $\pi_{n(n+1)}$ to a rescaled copy of $\bar Q_0$ in $\mathcal P_{0n}$. Thus, we can find
a $C_{n}'$-Lipschitz embedding of $\bar Q_{n+1}$ into $\bar Q_n$ that agrees with $\pi_{n(n+1)}$ on 
the outer squares of the rescaled copies of $\bar Q_1$ in the partition $\mathcal P_{1(n+1)}$,
and with $C_n'$ arbitrarily close to 1. 
Post-composing this map with $L_n$ we get a $C_{n+1}$-Lipschitz embedding $L_{n+1}$ of $\bar Q_{n+1}$ into $\bar Q_0$, where $C_{n+1}=C_nC_n'$. The sequence $(C_n)$ is monotone increasing and it can be chosen to converge to a constant $C>1$. 
By the Arzel\`a--Ascoli theorem the sequence of maps $(L_n\circ\pi_n)$ subconverges to a $C$-Lipschitz map $L_\infty$ of $S_2$ into $\bar Q_0\subset\R^2$. 

If $p$ and $q$ are two points in $S_2$ such that $\pi(p)\neq\pi(q)$, then from the definition of $L_\infty$ we clearly obtain that $L_\infty(p)\neq L_\infty(q)$.
The above construction also shows that if $J$ is a slit in $Q_m$ and $n\geq m$, then $L_n$ and $L_m\circ \pi_{mn}$ agree on $\pi_{mn}^{-1}(J)$, and thus $L_\infty$ is an embedding when restricted to the slit $\pi_m^{-1}(J)$. Therefore, $L_\infty$ is an embedding and let $S$ denote the image of $S_2$ under $L_\infty$.  For a slit $J$ in $S_2$, let $D_J$ denote the Jordan domain bounded by $L_\infty(J)$. Then $S=\bar Q_0\setminus \cup _{J}D_J$, where the union is over all slits of $S_2$.  Since $L_\infty$ is an embedding, the boundaries $\dee D_J$ are pairwise disjoint Jordan curves, disjoint from $\dee Q_0$. The sequence of diameters $({\rm diam}(\dee D_J))$ goes to 0 because $L_\infty$ is Lipschitz. The set $S$ has no interior because the slits are dense in $S_2$.
\qed

Let $L, R, T, B$ denote the \emph{left}, \emph{right}, \emph{top} and \emph{bottom} sides of the outer square of $S_2$, respectively, i.e., 
$$
\begin{aligned}
L&=\pi^{-1}(\{(0,y)\:0\leq y\leq1\}),\
R=\pi^{-1}(\{(1,y)\:0\leq y\leq1\}),\\ 
T&=\pi^{-1}(\{(x,1)\:0\leq x\leq1\}),\ 
B=\pi^{-1}(\{(x,0)\:0\leq x\leq1\}).
\end{aligned}
$$
Let $DS_2$ denote the \emph{double} of $S_2$ across the outer square, i.e., $DS_2$ is obtained by gluing two copies of $S_2$ along the sides of the outer square using isometries to identify the left, right, top, and bottom sides.  We refer to these two copies as the \emph{front} and \emph{back} copies.
As a double, the space $DS_2$ with the topology induced from $S_2$ is a Sierpi\'nski carpet as well, and that can be seen, for example, from Whyburn's characterization. \emph{Slits} of $DS_2$ are the slits of the two copies of $S_2$, and they form the family of all peripheral circles of $DS_2$.
The path metric $d_{S_2}$ on $S_2$ induces a  path metric $d_{DS_2}$ on $DS_2$.   Abusing notation slightly, we also denote by $\pi$ the projection of $DS_2$ onto $\bar Q_0$ induced by the projection $\pi$ of $S_2$ onto $\bar Q_0$: if $p\in DS_2$, then $p$ belongs to a copy of $S_2$ in $DS_2$ and $\pi(p)$ is the image of $p$ under $\pi\: S_2\to\bar Q_0$. From the definition of $DS_2$ we see that if $p$ belongs to both copies of $S_2$ in $DS_2$, the projections $\pi$ of each of these copies agree at $p$.

If $(X,d)$ is a metric space, $p\in X$, and $r> 0$, we denote by $B(p,r)$ the open ball in $X$ centered at $p$ with radius $r$, i.e.,
$$
B(p,r)=\{q\in X\: d(q,p)<r\}.
$$

\begin{lemma}\label{L:Incl}
 There exists a constant $c>0$ such that for every $p\in S_2$ and every $0<r\leq{\rm diam}(S_2)$, there exists $q\in\pi(S_2)$ with
$$
B(q,c\cdot r)\subseteq \pi(B(p,r))\subseteq B(\pi(p),r).
$$
\end{lemma}
\no
\emph{Proof.}
The right inclusion is obvious since
$\pi$ is 1-Lipschitz.
To show the left inclusion we first observe that since 
${\rm diam}(S_2)\leq 3$, for every $p\in S_2$ there exists $q=(1/2,1/2)\in \pi(S_2)$ such that 
\begin{equation}\label{E:Firststep}
B(q, {\rm diam}(S_2)/6)\subseteq\pi(B(p,{\diam}(S_2))).
\end{equation}
If $r$ is arbitrary, $0<r<{\rm diam}(S_2)$, there exists $n\in\N$ such that 
\begin{equation}\label{E:Secondstep}
{\rm diam}(S_2)\leq 2^nr<2\,{\rm diam}(S_2).
\end{equation}
The point $p$ belongs to a copy $S\subset S_2$ of $S_2$ rescaled by $1/2^n$ and thus we can apply~(\ref{E:Firststep}) to this copy to conclude that there exists $q\in \pi(S)$ such that 
$$
B(q,{\rm diam}(S)/6)\subseteq \pi(B(p,{\rm diam}(S))\cap S). 
$$
Since  ${\rm diam}(S)={\rm diam}(S_2)/2^n$, from~(\ref{E:Secondstep}) we conclude that
$$
B(q,r/12)\subseteq\pi(B(p,r)\cap S).
$$
Finally, since $B(p,r)\cap S\subseteq B(p,r)$, we obtain the left inclusion in the statement of the lemma with $c=1/12$.
\qed

The next lemma combined with the fact that every $\pi_n$ is Lipschitz implies that $\pi_n\: S_2\to \bar Q_n$ is a regular mapping, see~\cite[Definition~12.1]{DS97}.
\begin{lemma}\label{L:Reg}
There exists $C\geq 1$ such that for every $n\in\N\cup\{0\}$, for every $p\in S_2$, and $r>0$, the preimage $\pi_n^{-1}(B(\pi_n(p),r))$ can be covered by at most $C$ balls in $S_2$ of radii at most $C\cdot r$.
\end{lemma}
\no
\emph{Proof.}
Since $\pi=\pi_0$ factors as $\pi=\pi_{0n}\circ\pi_n$ and the maps $\pi_{0n}$ are 1-Lipschitz, it is enough to prove the lemma for $n=0$. Also it is enough to consider $r\leq1$. 

From compactness of $S_2$ it follows that there exists $C\geq 1$ such that if $1/4\leq r\leq 1$ and $p\in S_2$ is arbitrary, then the closure of $\pi^{-1}(B(\pi(p),r))$ can be covered by at most $C$ balls of radii at most one. 
Let $r$ be arbitrary now, $0<r<1/4$. There exists $n\in\N$ such that 
$$
1/4\leq2^nr<1/2. 
$$
The preimage $\pi^{-1}(B(\pi(p),r))$ is contained in the union of at most four rescaled by $1/2^n$ copies of $S_2$. Let $S$ be one of these copies. The intersection of $B(\pi(p),r)$ with $\pi(S)$ is contained in a ball $B$ centered at a point in $\pi(S)$ and whose radius is $r$. By the above, the preimage $\pi^{-1}(B)$ can be covered by at most $C$ balls in $S$ with radii at most $C\cdot r$. Since this holds for each of the four copies, the lemma follows.
\qed

The fact that each $\pi_n$ is 1-Lipschitz and Lemma~\ref{L:Reg} imply that there exists a constant $C\geq1$, independent of $n$, such that for any Borel set $E$ in $S_2$ we have
\begin{equation}\label{E:Compmeas}
\frac1{C}\mH^2(\pi_n(E))\leq\mH^2(E)\leq C\cdot\mH^2(\pi_n(E)),
\end{equation}
see~\cite[Lemma~12.3]{DS97}. 

If $B$ is a ball with center $p$ and radius $r$, and $\lambda$ is a positive constant, we denote by $\lambda B$ the ball centered at $p$ whose radius is $\lambda\cdot r$.
A metric measure space $(X,d,\mu)$ is said to be \emph{Ahlfors $Q$-regular} if there exists a constant $C\geq1$ such that
$$
\frac{r^Q}{C}\leq\mu(B(p,r))\leq C\cdot r^Q
$$
for all $p\in X$ and $0<r\leq{\rm diam}(X)$.
We say that $(X,d,\mu)$ is \emph{doubling} if there is a constant $C\geq1$ such that for every ball $B$ in $X$ we have
$$
\mu(2B)\leq C\cdot \mu(B).
$$
A metric Sierpi\'nski carpet $S$ is called \emph{porous} if there exists $C\geq1$ such that for every $p\in S$ and $0<r\leq{\rm diam}(S)$, there exists a peripheral circle $J$ in $S$ with $J\cap B(p,r)\neq\emptyset$ and
$$
\frac{r}{C}\leq {\rm diam}(J)\leq C\cdot r.
$$

\begin{proposition}\label{P:Basicprops}
 The Sierpi\'nski carpets $S_2$ and $DS_2$ with the path metrics and the Hausdorff 2-measures $\mH^2$ are compact, path-connected,  porous, Ahlfors 2-regular spaces, and in particular they are doubling.
\end{proposition}
\no
\emph{Proof.}
We only need to check the porosity and Ahlfors regularity for $S_2$. The other properties have already been established for $S_2$ and all of them extend to $DS_2$ in a straightforward way.

From the construction of $S_2$ it is clear that the projection $\pi(p)$ of every point $p$ in $S_2$ is within the Euclidean distance $\sqrt{2}/2^n$ from the projection of a slit of diameter $1/2^n,\ n\in\N$. Let $0< r\leq{\rm diam}(S_2)$ be arbitrary and $n\in\N$ be such that 
$$
\frac{\sqrt{2}}{2^n}<c\cdot r\leq\frac{2\sqrt{2}}{2^n}, 
$$
where $c$ is the constant from Lemma~\ref{L:Incl}.
The desired peripheral circle $J$ in the definition of porosity  is a peripheral circle of diameter $1/2^n$ whose projection is within distance $\sqrt{2}/2^n$ from the point $q$ as in Lemma~\ref{L:Incl}. This follows from the left inclusion of Lemma~\ref{L:Incl}.

For Ahlfors regularity, let $B(p,r)$ be any ball in $S_2$ with $0<r\leq {\rm diam}(S_2)$. From~(\ref{E:Compmeas}) we have
$$
a\cdot \mH^2(\pi(B(p,r)))\leq \mH^2(B(p,r))\leq b\cdot\mH^2(\pi(B(p,r))),
$$
where $a,b>0$ are independent of $p$ and $r$. Lemma~\ref{L:Incl} now gives a constant $C\geq1$, independent of $p$ and $r$, so that the left-hand side is at least $r^2/C$ and the right-hand side is at most $C\cdot r^2$.
\qed

\section{Non-vertical curve families}\label{S:Nv}

\no
We say that a curve $\gamma$ \emph{connects} two connected sets $E$ and $F$ in $X$ if $E\cup F\cup \overline\gamma$ is connected, where $\overline\gamma$ denotes the closure of the image of $\gamma$ in $X$. 
A curve $\gamma$ in $S_2$ or $DS_2$ is called \emph{vertical} if the first coordinate of $\pi(\gamma)$ is constant. Otherwise it is said to be \emph{non-vertical}.

If $\Gamma$ is a family of curves in a metric measure space $(X,d,\mu)$ and $Q>1$, its $Q$-\emph{modulus} is defined as
$$
{\rm mod}_Q(\Gamma)=\inf\bigg\{\int_X\rho^Qd\mu\bigg\},
$$
where the infimum is over all non-negative Borel functions $\rho$ on $X$ with
$$
\int_\gamma\rho\, ds\geq1
$$
for every locally rectifiable $\gamma\in\Gamma$. Here $ds$ denotes the arc-length element. Such a  function $\rho$  is  referred to as a \emph{mass distribution} on $X$ and the quantity $\int_X\rho^Qd\mu$ is called the \emph{total mass} of $\rho$. 

It is known that the $Q$-modulus is 
\emph{monotone}, i.e., 
 $$
 {\rm mod}_Q(\Gamma_1)\leq {\rm mod}_Q(\Gamma_2)
 $$
 for $\Gamma_1\subseteq \Gamma_2$, and
\emph{subadditive}, i.e., if $\Gamma=\cup\Gamma_k$, then
$$
{\rm mod}_Q(\Gamma)\leq\sum{\rm mod}_Q(\Gamma_k).
$$
If $Q=2$, we write ${\rm mod}$ instead of ${\rm mod}_2$. 

Recall that a homeomorphism $f$ from a planar domain $D$ onto a planar domain $\tilde D$ is said to be $K$-\emph{quasiconformal}, $K\geq 1$, if it is ACL (absolutely continuous on almost every line) and
$$
\bigg|\frac{\dee f}{\dee \bar z}\bigg|\leq k\bigg|\frac{\dee f}{\dee z}\bigg|
$$
almost everywhere in $D$, where $k=(K-1)/(K+1)$. 
If $K=1$, a map is called \emph{conformal}. See~\cite{lA66} for background on quasiconformal maps. 
If $f\: D\to\tilde D$ is a $K$-quasiconformal map, $\Gamma$ is a curve family in $D$, and $\tilde \Gamma=f(\Gamma)$, then there is a constant $C\geq1$ that depends only on $K$ such that
\begin{equation}\label{E:Qimod}
\frac1{C}{\rm mod}(\Gamma)\leq{\rm mod}(\tilde\Gamma)\leq C\cdot {\rm mod}(\Gamma).
\end{equation}
If $f$ is conformal, then ${\rm mod}(\tilde\Gamma)={\rm mod}(\Gamma)$.
The quasi-invariance of the modulus as in~(\ref{E:Qimod}) has been greatly extended to cover quasiconformal or quasisymmetric maps between more general metric measure spaces, see, e.g.,~\cite{jH01}. In particular, if $X$ and $Y$ are locally compact, connected, Ahlfors $Q$-regular metric spaces  and $f\: X\to Y$ is a quasisymmetric map, then there exists a constant $C\geq1$ such that the inequalities as in~(\ref{E:Qimod}), with ${\rm mod}$ replaced by ${\rm mod}_Q$, hold for every curve family $\Gamma$ in $X$,  see~\cite{jT98}.

\begin{lemma}\label{L:Horfam}
For every $\epsilon>0$ there exists $n\in\N$ and a conformal map $\Phi_n$ from $Q_n$ onto a multiply connected domain $D_n$ in an open rectangle  $(0,M_n)\times(0,1)$ such that the homeomorphic extension of $\Phi_n$ takes the vertices $(0,0), (0,1), (1,1), (1,0)$ to $(0,1), (0,1), (M_n,1), (M_n,0)$, respectively, and
$$
\frac1{M_n}<\epsilon.
$$
\end{lemma}
\no
\emph{Proof.}
Let $M>0$ and let $\phi_M$ be the unique conformal map from $(0,M)\times(0,1)$ to $(0,\tilde M)\times(0,1)$ so that the points $(0,0)$, $(0,1)$, $(M,1/2)$, $(M,0)$ on the boundary go under the homeomorphic extension, also denoted by $\phi_M$, to the vertices $(0,0)$, $(0,1)$, $(\tilde M, 1)$, $(\tilde M,0)$, respectively. From monotonicity of the modulus we obtain $\tilde M\geq M$.

Now we define a conformal map $\Phi_n$ of $Q_n$ onto a multiply connected domain in $(0,M_n)\times(0,1)$ inductively by $\Phi_0=\id, M_0=1$, and 
$$
\Phi_{n+1}(p)=\frac12\phi_{M_n}(\Phi_{n}(2p)),\quad {\rm for}\ p\in\frac12 Q_n\subset Q_{n+1},
$$
extended by reflections and homeomorphically to the rest of $Q_{n+1}$.
Observe that for every $n\in \N\cup\{0\}$, the homeomorphic extension of $\Phi_n$ takes the point $(1,1/2)$ to $(M_n,1/2)$, and 
$$
M_{n+1}=\phi_{M_n}(M_n).
$$
Monotonicity of the modulus gives that the sequence $(M_n)$ is non-decreasing and we assume for contradiction that $M_\infty=\sup\{M_n\}<\infty$. 

The horizontal stretching 
$
T_n\: (0,M_n)\times(0,1)\to(0,M_\infty)\times(0,1)
$
is a $(M_\infty/M_n)$-quasiconformal map.
Consider the map 
$$
T_{n+1}\circ\phi_{M_n}\circ T_n^{-1}\: (0,M_\infty)\times(0,1)\to(0,M_\infty)\times(0,1).
$$
It is $(M_\infty^2/(M_nM_{n+1}))$-quasiconformal, its homeomorphic extension fixes $(0,0), (0,1), (M_\infty,0)$,  and takes $(M_\infty, 1/2)$ to $(M_\infty,1)$. By choosing $n$ large enough, we obtain a contradiction as follows. Conjugating each map $T_{n+1}\circ\phi_{M_n}\circ T_n^{-1}$ by the same conformal map $\phi$ of $(0,M_\infty)\times(0,1)$ onto the unit disc, we obtain a sequence $(f_n)$ of self-maps of the unit disc such that $f_n$ is $(M_\infty^2/(M_nM_{n+1}))$-quasiconformal and the continuous extension of $f_n$ to the boundary fixes three distinct points, say $1, i, -1$. The map $\phi$ has a continuous extension to the boundary and takes $(M_\infty, 1/2)$ to $\xi$ and $(M_\infty,1)$ to $\zeta$, where $\xi$ and $\zeta$ are two distinct points on the boundary of the unit disc. Using the Schwarz reflection principle we get a sequence $(\tilde f_n)$ of self-maps of the Riemann sphere, where $\tilde f_n$ is $(M_\infty^2/(M_nM_{n+1}))$-quasiconformal and it fixes $1, i, -1$. This is a compact family and  
thus it has a subsequence that converges uniformly to a $(M_\infty^2/(M_kM_{k+1}))$-quasiconformal map $f$ for every $k$. Therefore, $f$ is 1-quasiconformal, and hence a M\"obius transformation. Since it fixes three distinct points, $f$ is the identity. However, all of the maps $\tilde f_n$, and hence $f$, take $\xi$ to $\zeta$, a contradiction. 
\qed

\begin{corollary}\label{C:Nv}
Let $\Gamma$ be a family of non-vertical curves in $(S_2, d_{S_2},\mH^2)$ or in $(DS_2, d_{DS_2},\mH^2)$. Then
$$
{\rm mod}(\Gamma)=0.
$$ 
\end{corollary}
\no
\emph{Proof.}
We can write 
$
\Gamma=\cup\Gamma_k
$,
where $\Gamma_k$ consists of all curves $\gamma$ in $\Gamma$ such that the oscillation of the first coordinate of $\pi(\gamma)$, i.e., the difference between the supremum and the infimum, is at least $1/k$. 
Since the modulus is subadditive, it is enough to show that ${\rm mod}(\Gamma_k)=0$ for each $k\in\N$.

Let $\Gamma$ be a curve family in $S_2$ and let $\epsilon>0$ be arbitrary.
We choose $m\in\N$ such that $1/2^{m}<1/(2k)$. Then each curve in $\Gamma_k$ connects the two complementary components of 
$$
S_2\cap\pi^{-1}\bigg(\bigg\{\frac{l}{2^m}< x<\frac{l+1}{2^m}\bigg\}\cap \bar Q_0\bigg)
$$
in $S_2$ for some $0\leq l\leq 2^m-1$.
By Lemma~\ref{L:Horfam}, we can find $n\in\N$ such that there exists a conformal map $\Phi_n$
of $Q_n$ onto a multiply connected domain in $(0,M_n)\times(0,1)$ whose homeomorphic  extension to the boundary takes $(0,0), (0,1), (1,1), (1,0)$ to $(0,1), (0,1), (M_n,1), (M_n,0)$, respectively, and $2^{2m}/{M_n}<\epsilon$.
The map
$$
p\mapsto \frac1{2^m}\Phi_n(2^mp),\quad p\in \frac1{2^m} Q_n,
$$
extends by reflections and homeomorphically to a conformal map $\Phi_{m,n}$ from $Q_N$ onto a multiply connected domain in the rectangle $(0,M_n)\times(0,1)$, where $N=m+n$. 
Notice that if $\gamma$ is in $\Gamma_k$, the oscillation of the first coordinate of $\pi_N(\gamma)\cap Q_N$ is at least $1/k$, and thus the total length of $\Phi_{m,n}(\pi_N(\gamma)\cap Q_N)$, i.e., the sum of the lengths of its connected components, is at least $M_n/2^m$. The map $\tilde\Phi_{m,n}$ defined by
$$
p\mapsto\frac{2^m}{M_n}\Phi_{m,n}(p),\quad p\in Q_N,
$$
is conformal from $Q_N$ onto a multiply connected domain $D_{\epsilon,k}$ in a rectangle $(0,w_k)\times(0,h_{\epsilon, k})$, where $w_k=2^m$ and $h_{\epsilon,k}=2^m/M_n$. If $\gamma\in\Gamma_k$, the total length of $\tilde\Phi_{m,n}(\pi_N(\gamma)\cap Q_N)$ is at least one. 

Let $\rho_{N}$ be the mass distribution on $\bar Q_N$ given by $|\tilde\Phi'_{m,n}|$ in $Q_N$ and 0 on the slits of $\bar Q_N$. The mass distribution $\rho_{N}$ is admissible for $\pi_N(\Gamma_k)$ because for every $\gamma\in\Gamma_k$
$$
\int_{\pi_N(\gamma)}\rho_{N}ds=\int_{\pi_N(\gamma)\cap Q_N}|\tilde\Phi'_{m,n}|ds=\int_{\tilde\Phi_{m,n}(\pi_N(\gamma)\cap Q_N)}ds,
$$
and the last integral is the total length of $\tilde\Phi_{m,n}(\pi_N(\gamma)\cap Q_N)$, which is at least one. The total mass of $\rho_{N}$ is 
$$
\begin{aligned}
\int_{\bar Q_N}\rho_{N}^2d\mH^2&=\int_{Q_N}|\tilde\Phi'_{m,n}|^2 dxdy\\
&\leq\int_{(0,w_k)\times(0,h_{\epsilon, k})}dudv= w_k\cdot h_{\epsilon,k}=2^{2m}/M_n<\epsilon.
\end{aligned}
$$
In particular,
${\rm mod}(\pi_N(\Gamma_k))<\epsilon$. 

Now consider $\rho(p)=\rho_{N}(\pi_N(p)),\ p\in S_2$, a mass distribution on $S_2$. Since $\pi_N$ is 1-Lipschitz, we obtain
$$
\int_\gamma\rho ds\geq\int_{\pi_N(\gamma)}\rho_Nds\geq1,\quad\forall\gamma\in\Gamma_k.
$$
By~(\ref{E:Compmeas}) there exists a constant $C\geq1$, independent of $N$, such that for any Borel set $E$ in $S_2$ we have
$$
\frac1{C}\mH^2(\pi_N(E))\leq\mH^2(E)\leq C\cdot\mH^2(\pi_N(E)).
$$
Therefore,
$$
\int_{S_2}\rho^2d\mH^2\leq C\int_{\bar Q_N}\rho_N^2d\mH^2<C\cdot\epsilon.
$$
Since $C$ is independent of $N$, the conclusion of the lemma follows for the space $S_2$.
The case of $DS_2$ requires only minor modifications.
\qed

\section{Ahlfors regularity of the image}\label{S:Ahlreg}

\no
For a subset $E$ in $X$ we denote by $\chi_E$ its characteristic function. 
In what follows we need the following lemma, see~\cite[Exercise 2.10]{jH01}.
\begin{lemma}\label{L:Boj}
Let $(B_i)$ be a countable collection of pairwise disjoint balls in a doubling metric measure space $(X,d,\mu)$, let $(a_i)$ be non-negative numbers, and let $\lambda\geq1$ be arbitrary. Then
$$
\int_X\bigg(\sum_ia_i\chi_{\lambda B_i}\bigg)^2d\mu\leq C\int_X\bigg(\sum_ia_i\chi_{B_i}\bigg)^2d\mu,
$$
where $C$ depends only on $\lambda$ and the doubling constant of $\mu$.
\end{lemma}
\no
\emph{Proof.}
It follows easily from the doubling property of the measure $\mu$ that for every $\lambda\geq0$ there exists a constant $C_\lambda$ that depends only on $\lambda$ and the doubling constant of $\mu$, such that
$$
\mu(\lambda B)\leq C_\lambda\mu(B)
$$
for every ball $B$. Because of this, in what follows we may assume that $\mu(B_i)>0$ for all $i$.

Let $\phi\in L^2=L^2(X,\mu)$. We denote the non-centered maximal function of $\phi$ by $M_{nc}(\phi)$, i.e., 
$$
M_{nc}(\phi)(p)=\sup_{B}\bigg\{\frac1{\mu(B)}\int_{B}|\phi|d\mu\bigg\},
$$
where the supremum is over all open balls $B$  containing $p$.
Then 
$$
\begin{aligned}
\bigg|\int_X\sum_{i} a_i\chi_{\lambda B_i}\phi(q)d\mu(q)\bigg|
&=\bigg|\sum_{i}a_i\int_{\lambda B_i}\phi(q)d\mu(q)\bigg|\\
&\leq\sum_{i}a_i\frac{\mu(\lambda B_i)}{\mu(B_i)}\int_{B_i}M_{nc}(\phi)(p)d\mu(p)\\
&\leq C_\lambda\int_X\sum_{i} a_i\chi_{B_i}M_{nc}(\phi)(p)d\mu(p)\\
&\leq C_\lambda\bigg|\bigg|\sum_{i}a_i\chi_{B_i}\bigg|\bigg|_{L^2}\cdot||M_{nc}(\phi)||_{L^2}.
\end{aligned}
$$
Since $\mu$ is doubling, we have 
$$
M_{nc}(\phi)\leq C_\mu M(\phi),
$$
where $C_\mu $ is the doubling constant for $\mu$ and $M$ is the maximal function defined by
$$
M(\phi)(p)=\sup_{r>0}\bigg\{\frac1{\mu(B(p,r))}\int_{B(p,r)}|\phi|d\mu\bigg\}.
$$
As is well-known, see, e.g., \cite[Theorem 2.2]{jH01},  the maximal function satisfies the inequality
$$
||M(\phi)||_{L^2}\leq C_2||\phi||_{L^2}.
$$
Combining the above estimates and using the duality of $L^2$ we obtain the desired inequality with $C=C_\lambda^2 C_\mu^2C_2^2$.
\qed

\begin{lemma}\label{L:ModS2}
Let $\Gamma$ be the family of all vertical curves in $S_2$ that connect $T$ and $B$. Then
$$
0<{\rm mod}(\Gamma)<\infty.
$$
\end{lemma}
\no
\emph{Proof.} It is clear that ${\rm mod}(\Gamma)$ is finite because $\rho\equiv 1$ is an admissible mass distribution for $\Gamma$.
Let $\rho$ be an arbitrary admissible mass distribution for $\Gamma$. 
For every $x\in B$ let $\gamma_x$ denote a vertical curve in $\Gamma$ that contains $x$. Let $B'$ be the subset of $B$ that consists of all $x\in B$ so that  $\gamma_x$ does not intersect any slits of $S_2$. Slightly abusing notation we let $B'\times L$ denote the measurable subset of $S_2$ that consists of all points in $S_2$ that belong to $\gamma_x$ for some $x\in B'$. 
There is a natural product measure $\mu$ on $B'\times L$, which comes from the Hausdorff 1-measures on $B'$ and $L$. It follows easily from~(\ref{E:Compmeas}) that $\mu$ is roughly comparable to $\mH^2$ on $B'\times L$, i.e., there exists  a constant $C\geq 1$ such that for any Borel set $E$ in $B'\times L$ we have
\begin{equation}\label{E:Productmeas}
\frac1{C}\mH^2(E)\leq\mu(E)\leq C\cdot\mH^2(E).
\end{equation}

Now, for any $x\in B'$, 
$$
\int_{\gamma_x}\rho\, ds\geq 1.
$$
Using the Cauchy--Schwarz inequality and the fact that the length of each $\gamma_x$ is one, we obtain
$$
\int_{\gamma_x}\rho^2ds\geq1.
$$
Integrating over $x\in B'$ gives
$$
\int_{B'}\bigg(\int_{\gamma_x}\rho^2ds\bigg)dx\geq1.
$$
Applying Fubini's theorem we obtain
$$
\int_{B'\times L}\rho^2d\mu\geq1.
$$ 
This inequality combined with~(\ref{E:Productmeas}) now gives
$$
\int_{S_2}\rho^2d\mH^2\geq\frac1{C}.
$$
Since $\rho$ is an arbitrary admissible mass distribution, we conclude that ${\rm mod}(\Gamma)\geq1/C$.
\qed

\begin{lemma}\label{L:Ahl}
Let $f\: S_2\to S_2$, or $f\: DS_2\to DS_2$, be a quasisymmetric embedding. Then the image $f(S_2)$, or $f(DS_2)$, is Ahlfors 2-regular. 
\end{lemma}
\no
\emph{Proof.} As before, we only treat the case of $S_2$. A proof for $DS_2$ follows the same lines with minor modifications.

We need to show that there exists $C\geq 1$ such that
$$
r^2/C\leq\mH^2(B(\tilde p,r))\leq C\cdot r^2
$$
for all $\tilde p\in f(S_2)$ and $0<r\leq{\rm diam}(f(S_2))$, where $B(\tilde p,r)$ denotes the ball in $(f(S_2),d_{S_2})$. The upper bound follows immediately because $S_2$ is Ahlfors 2-regular. Assume for contradiction that there exists a sequence $(\tilde p_n),\ \tilde p_n\in f(S_2)$, and a sequence $(\tilde r_n),\ 0<\tilde r_n\leq{\rm diam}(f(S_2))$, such that for the ball $B(\tilde p_n,\tilde r_n)$ in $(f(S_2),d_{S_2})$ we have
\begin{equation}\label{E:Hm}
\frac{\mH^2(B(\tilde p_n, \tilde r_n))}{\tilde r_n^2}\to0\quad {\rm as}\ n\to\infty.
\end{equation}
Let $\tilde B_n$ denote the ball $B(\tilde p_n, \tilde r_n)$ in the metric space 
$$
\bigg(f(S_2), \frac1{\tilde r_n}d_{S_2}\bigg).
$$
This ball has radius one and the limit in~(\ref{E:Hm}) gives
$$
\mH^2(\tilde B_n)\to0 \quad {\rm as}\ n\to\infty.
$$
Let $B_n$ denote the set $f^{-1}(B(\tilde p_n,\tilde r_n))$ in the metric space 
$$
\bigg(S_2,\frac1{r_n}d_{S_2}\bigg),
$$
where $r_n>0$ is chosen so that the diameter of $B_n$ is one. We denote by $f_n$ the map $f$ between the rescaled spaces.

Since $f$ is quasisymmetric, say with a distortion function $\eta$, and pre- or post-compositions with scalings do not change the distortion function, $f_n$ is $\eta$-quasisymmetric. Thus, there is a constant $M>1$, independent of $n$, such that
$$
B(p_n,1/{M})\subseteq B_n\subseteq B(p_n,M),
$$ 
where $p_n=f^{-1}(\tilde p_n)$. 

By Lemma~\ref{L:Incl}, there exist $c>0$ and $q_n\in\pi(S_2)$ such that 
\begin{equation}\label{E:Cont}
B(q_n, {c}/{M})\subseteq\pi(B(p_n,{1}/{M})).
\end{equation}
Let $S$ be a copy of $S_2$ contained in $B(p_n,1/M)$ and having the maximal diameter. By~(\ref{E:Cont}),  the side length of the outer square of $S$ is at least $\delta>0$ that is independent of $n$. 
Let $\Gamma_n$ be the family of all vertical curves in $S$ that connect the top and the bottom sides of the outer square of $S$. The diameter of each $\gamma$ in $\Gamma_n$ is 
at least 
$\delta$. In addition, by Lemma~\ref{L:ModS2} and by the obvious fact that the modulus is invariant under scalings  we have ${\rm mod}(\Gamma_n)=\sigma>0$, where $\sigma$ is independent of $n$.

Let $\tilde\Gamma_n=f_n(\Gamma_n)$. 
Since all $f_{n}$ are $\eta$-quasisymmetric with the same $\eta$, the diameter of $B_n$ is one and the radius of $\tilde B_n$ is one, Proposition~10.8 in~\cite{jH01} implies that there exists $\tilde\delta>0$, independent of $n$, such that for every  $\gamma\in\Gamma_n$ the image $\tilde\gamma=f_n(\gamma)$ has the diameter at least $\tilde\delta$. 

The proof now follows the lines of that of Theorem~15.10 in~\cite{jH01}. Let $\epsilon>0$ and let $n$ be chosen so that $\mH^2(\tilde B_n)<\epsilon$. We fix a disjoint collection of balls $(\tilde B_i'),\ \tilde B_i'=B(\tilde p_i', \tilde s_i')$, in $\tilde B_n$ such that
the collection $(5\tilde B_i')$ covers $\tilde B_n$ and 
$$
\sum_i(\tilde s_i')^2<\epsilon.
$$
This is possible by Theorem~1.2 in ~\cite{jH01}. Since $f_n$ is $\eta$-quasisymmetric, there exists $H\geq1$, independent of $n$ and $i$, and a collection of balls $(B_i'),\ B_i'=B(p_i', s_i')$, such that 
$$
B_i'\subseteq f_n^{-1}(\tilde B_i')\subseteq f_n^{-1}(5\tilde B_i')\subseteq HB_i'.
$$ 
By choosing $\epsilon$ small enough, we may and will assume that $4H\cdot s_i'<\delta$.

Now we consider a mass distribution $\rho$ on $B_n$ defined by
$$
\rho=\sum_i \frac{10\tilde s_i'}{H\tilde\delta s_i'}\chi_{2HB_i'}.
$$
Let $\gamma\in\Gamma_n$ and let $i$ be an index such that 
$$
f_n(\gamma)\cap (5\tilde B_i')\neq\emptyset.
$$
Then $\gamma\cap HB_i'\neq\emptyset$, and since 
$$
{\rm diam}(\gamma)\geq\delta>4H\cdot s_i'\geq{\rm diam}(2HB_i'),
$$ 
the curve $\gamma$ cannot be completely contained in $2HB_i'$. This gives 
$$
{\rm length}(\gamma\cap2HB_i')\geq H\cdot s_i'.
$$
Therefore,
$$
\begin{aligned}
\int_\gamma\rho\, ds&=\sum_i\frac{10\tilde s_i'}{H\tilde\delta s_i'}{\rm length}(\gamma\cap2HB_i')\\
&\geq\frac{1}{\tilde\delta}\sum_{i\: f_n(\gamma)\cap (5\tilde B_i')\neq\emptyset}10\tilde s_i'\\
&\geq \frac1{\tilde\delta}\, {\rm diam}(f_n(\gamma))
\geq1.
\end{aligned}
$$

Now, using Lemma~\ref{L:Boj}, we obtain
$$
\begin{aligned}
{\rm mod}(\Gamma_n)&\leq\int_{B_n}\rho^2d\mH^2\\
&=\int_{B_n}\bigg(\sum_i\frac{10\tilde s_i'}{H\tilde\delta s_i'}\chi_{2HB_i'}\bigg)^2d\mH^2\\
&\leq C\int_{B_n}\bigg(\sum_i\frac{10\tilde s_i'}{H\tilde\delta s_i'}\chi_{B_i'}\bigg)^2d\mH^2\\
&\leq C\sum_i\frac{(10\tilde s_i')^2}{(H\tilde \delta s_i')^2}(s_i')^2\\
&=\frac{100\cdot C}{(H\tilde\delta)^2}\sum_i(\tilde s_i')^2
<\frac{100\cdot C}{(H\tilde\delta)^2}\epsilon.
\end{aligned}
$$
This contradicts the fact that ${\rm mod}(\Gamma_n)=\sigma$.
\qed

\begin{corollary}\label{C:Verttovert}
Every quasisymmetric embedding $f$ of $S_2$ or $DS_2$ into itself takes every vertical curve to a vertical curve.
\end{corollary}
\no
\emph{Proof.}
Let $\Gamma_{v\to nv}$  be the curve family in $S_2$ or $DS_2$ that consists of vertical curves mapped by $f$ to non-vertical curves. If $\Gamma_{v\to nv}$ were a non-empty family, then
Corollary~\ref{C:Nv} would imply that ${\rm mod}(f(\Gamma_{v\to nv}))=0$. By Lemma~\ref{L:Ahl}, $f(S_2)$, or $f(DS_2)$, is Ahlfors 2-regular, and therefore $f$ quasi-preserves the 2-modulus as defined in~(\ref{E:Qimod}), see~\cite{jT98}. We conclude that the 2-modulus of $\Gamma_{v\to nv}$ is zero as well. Combined with the proof of Lemma~\ref{L:ModS2}, this implies that for almost every $x$ in $B\subset S_2\subset DS_2$, any vertical curve that contains $x$ is mapped to a vertical curve.
This readily implies that every vertical curve in $S_2$ or $DS_2$ is mapped by $f$ to a vertical curve, i.e., $\Gamma_{v\to nv}$ is, in fact, empty.
\qed

\section{Co-Hopfian property}

\no
Here we finish the proof of Theorem~\ref{T:QsCoHopf}, and hence Theorem~\ref{T:QiCoHopf}, by proving the following theorem. A vertical curve in $DS_2$ is called \emph{closed} if it is homeomorphic to a circle.
\begin{theorem}\label{T:CohopF}
The Sierpi\'nski carpet $DS_2$ with the path metric is quasisymmetrically co-Hopfian.
\end{theorem}
\no
\emph{Proof.}
Let $f$ be a quasisymmetric embedding of $DS_2$ into itself and let $\Gamma_v$ be the family of all closed vertical curves in $DS_2$.  By Corollary~\ref{C:Verttovert}, $f$ maps $\Gamma_v$ to a family of vertical curves, and they must be closed since $f$ is a continuous embedding. A closed vertical curve that intersects the slits of $DS_2$ of the largest diameter (which is 1/2) must be mapped by $f$ to a curve with the same property because these are the only closed vertical curves that intersect only two of the slits of $DS_2$. There are four such curves and $f$ permutes them. Likewise, $f$ permutes closed vertical curves that intersect slits of diameter 1/4 and so forth. Since slits are dense in $DS_2$, the map $f$ is onto.
\qed

\section{Quasisymmetry groups}\label{S:Qsgroups}

\no
In this section we give full descriptions of the quasisymmetry groups of $S_2$ and $DS_2$. 
This section was added in the revised version of the paper after the question of describing these groups had been raised by Mario~Bonk and  the anonymous referee.

As before, we denote the left, right, top, and bottom sides of the outer square of $S_2$ by $L, R, T$, and $B$, respectively. Slightly abusing the notation, we use the same letters to denote the corresponding subsets of the double $DS_2$. Every slit $s$ of $S_2$ or $DS_2$ has two distinguished points on it: 
one closest to $T$ and the other to $B$.
We refer to them as the \emph{top-most} and the \emph{bottom-most} points of the slit $s$.

\begin{theorem}\label{T:QsgpS2}
The group of quasisymmetric self-maps of $S_2$ is the group of isometries of $S_2$, that is the finite dihedral group ${\bf D}_2$ consisting of four elements (it is isomorphic to $\Z_2\times\Z_2$).
\end{theorem}
\no
\emph{Proof.}
Let $R_r$ denote the rotation of $S_2$ induced by the rotation of the plane by $180^\circ$, and let $R_v$ and $R_h$ denote the reflections of $S_2$ induced by the reflections of the plane with respect to  $\{x=1/2\}$ and $\{y=1/2\}$, respectively. It is clear that $R_r,R_v, R_h$ are isometries of $S_2$, that $R_h=R_r\circ R_v$, and that these elements generate the dihedral group ${\bf D}_2$. 
According to Corollary~\ref{C:Verttovert}, every quasisymmetry $f$ takes vertical curves to vertical curves. This readily implies that $f$ takes slits to slits and the outer square is fixed setwise. 
Moreover, if $\{v_1,v_2,v_3,v_4\}$ is the ordered set of the four vertices of the outer square starting from the top left and moving clockwise, then $f$ either leaves this set pointwise fixed or maps it to one of the following ordered configurations: $\{v_3,v_4,v_1,v_2\}, \{v_2,v_1,v_4,v_3\}, \{v_4,v_3,v_2,v_1\}$. By composing $f$ with one of $R_r, R_v, R_h$ if needed, it is enough to show  that if the set of vertices of the outer square is fixed pointwise by $f$, then $f$ is the identity. 

Assuming that $f$ fixes the set of vertices of the outer square pointwise, we conclude that $f$ preserves each of $L, R, T$, and $B$ setwise and is an orientation preserving homeomorphism when restricted to each of these sets. The proof of Theorem~\ref{T:CohopF} shows that every vertical curve that connects $T$ to $B$ and intersects slits of diamater $1/2^n$ gets mapped onto a curve of the same type. Moreover, since the restriction of $f$ onto $T$ and $B$ preserves the orientations, each such curve is mapped onto a vertical curve with the same endpoints on $T$ and $B$. This implies that the top-most and the bottom-most points of every slit are fixed points of $f$. From the construction of $S_2$ it follows that the set of such points is dense in $S_2$, and since $f$ is continuous, it must be the identity.
\qed

\begin{lemma}\label{L:IsomDS2}
The group of isometries of $DS_2$ is isomorphic to $(\Z_2)^3$.
\end{lemma}
\no
\emph{Proof.}
Let $f$ be an isometry of $DS_2$. According to Corollary~\ref{C:Verttovert}, $f$ maps every vertical curve to a vertical curve. The vertical curves $L$ and $R$ are distinguished in the sense that they are not contained in closed vertical curves. Thus, the union $L\cup R$ is setwise fixed by $f$ and the set of the endpoints of $L$ and $R$ is mapped to itself. Also, since $f$ is an isometry and $T$ and $B$ are clearly geodesics, we have that $f$ fixes $T\cup B$, and again the endpoints of $T$ and $B$ go to the endpoints. Therefore, as in the case of $S_2$ above, if $\{v_1, v_2, v_3, v_4\}$ is the ordered set of the four vertices of the outer square as in Theorem~\ref{T:QsgpS2}, the map $f$ either leaves this set pointwise fixed or maps it to one of the following configurations: $\{v_3,v_4,v_1,v_2\}, \{v_2,v_1,v_4,v_3\}, \{v_4,v_3,v_2,v_1\}$. The isometries $R_r, R_v, R_h$ of $S_2$ defined in Theorem~\ref{T:QsgpS2} extend in a natural way to isometries of $DS_2$, and we denote them also by $R_r, R_v, R_h$. By possibly composing $f$ with one of these maps, we may assume that $f$ fixes the set $\{v_1, v_2, v_3, v_4\}$ pointwise, and therefore fixes each of $L,R, T, B$ setwise.   

This leaves two possibilities: either each of the copies, the front and the back, of $S_2$ in $DS_2$ stays setwise fixed, or they interchange. There is an isometry $R_{fb}$ of $DS_2$ that fixes $L\cup R\cup T\cup B$ pointwise and interchanges these two copies. By composing the map $f$ with such an isometry if needed, we may assume that $f$ leaves each of the two copies of $S_2$ in $DS_2$ setwise fixed. Now we may apply  Theorem~\ref{T:QsgpS2} to $f$ restricted to the front and the back copies to conclude that $f$ must be the identity. Therefore, the group of isometries of $DS_2$ is generated by $R_r, R_v$, and $R_{fb}$. 
It follows immediately that this group is isomorphic to $\Z_2\times\Z_2\times\Z_2$.
\qed

To describe the group of quasisymmetric self-maps of $DS_2$ we need to introduce an auxiliary  group. Let $\mathcal L$ be the additive group of real-valued Lipschitz functions $h$ on $[0,1]$ such that $h(0)=0$, $h(1)$ is an integer, and if $k/2^n,\ n\in\N,  k=1,\dots, 2^n-1$, is a reduced dyadic number, then 
$$
h\bigg(\frac{k}{2^n}\bigg)=\frac{2l}{2^{n}}
$$
for some $l\in\Z$ depending on $n$ and $k$. Here a dyadic number $k/2^n$ is \emph{reduced} if $k$ is odd. It is immediate that $\mathcal L$ is an abelian group. This group is uncountable. Indeed, let $h_0$ be a real-valued function on $[0,1]$ defined by
$$
h_0(t)=\begin{cases}0,\quad {\rm for}\ 0\leq t\leq1/2\\ \frac12-2|t-\frac34|,\quad {\rm for}\ 1/2\leq t\leq 1.
\end{cases}
$$
It is clear that $h_0$ is 2-Lipschitz and an elementary calculation shows that $h_0$ is in $\mathcal L$. The map $h_0$ extends by 0 to a continuous function on the whole real line. Now define a family of functions $h_{\epsilon}$ by the formula
$$
h_\epsilon(t)=\sum_{m=0}^\infty \epsilon_m\frac1{2^m}h_0({2^m}t),\quad t\in[0,1],
$$
where each $\epsilon_m$ is equal to either 0 or 1. The support of $h_0(2^mt)$ is $[1/2^{m+1},1/2^m]$, and therefore the series converges to a 2-Lipschitz function with $h_\epsilon(0)=h_\epsilon(1)=0$. If $k/2^n$ is a reduced dyadic number in $[0,1]$, there exists a unique $m\in\N\cup\{0\}$ such that $k/2^n\in(1/2^{m+1},1/2^m]$. Therefore, 
$$
h_\epsilon(k/2^n)=\frac{\epsilon_mh_0\big(\frac{2^mk}{2^n}\big)}{2^m}=\frac{\epsilon_m\frac{2l}{2^{n-m}}}{2^m}=\frac{2l\epsilon_m}{2^n}
$$
for some integer $l$. Thus, each $h_\epsilon$ is an element of $\mathcal L$, and hence $\mathcal L$ is uncountable.

Each closed vertical curve $\gamma$  in $DS_2$ can be given an orientation so that when one moves along $\gamma$ on the front copy of $S_2$ in $DS_2$ the distance from $B$ increases, and on the back copy it decreases. Intuitively, when one moves along $\gamma$ in the direction it is oriented, $L$ always ``stays on the left". We say that a quasisymmetric map $f$ of $DS_2$ is a \emph{rotation} on a closed vertical curve $\gamma$ if it maps $\gamma$ 
locally 
isometrically onto a closed vertical curve preserving the orientation. Note that we do not require $f$ to preserve $\gamma$.

\begin{theorem}\label{T:QsgpDS2}
The group $\mathcal{QS}$ of quasisymmetric self-maps of $DS_2$ coincides with the group of bi-Lipschitz self-maps of $DS_2$ and is the semidirect product  $\mathcal{I}\ltimes \mathcal{V}$, where $\mathcal I$ and $\mathcal V$ are the following subgroups. The subgroup $\mathcal I$ is the group of isometries of $DS_2$ and therefore it is isomorphic to $(\Z_2)^3$. The subgroup $\mathcal V$ is the group of all quasisymmetric self-maps that are the identity on $L$, coincide with isometries when restricted to $R$, and are rotations on every closed vertical curve. The group $\mathcal V$ is isomorphic to $\mathcal L$, and thus it is abelian and uncountable. 
\end{theorem}
\no
\emph{Proof.}
Let $f$ be an arbitrary quasisymmetric self-map of $DS_2$. 
As in the proof of Lemma~\ref{L:IsomDS2},
Corollary~\ref{C:Verttovert} gives that $f$ takes every vertical curve to a vertical curve and, in particular, $L\cup R$ is setwise fixed. By post-composing $f$ with an isometry of $DS_2$ we may therefore assume that $L$ is setwise fixed and its endpoints are fixed points.

From continuity of $f$ it follows that $f$ either simultaneously preserves or simultaneously reverses  the orientations of closed vertical curves. In the latter case we can post-compose $f$ with $R_{fb}$ and therefore we may assume that $f$ restricted to every closed vertical curve is an orientation preserving homeomorphism. We will show that such $f$ must be an element of $\mathcal V$. 

The map $f$ cannot interchange the ``order" of closed vertical curves, namely if $\gamma_1$ and $\gamma_2$ are two vertical curves so that the distance from $\gamma_1$ to $L$ is less than the distance from $\gamma_2$ to $L$, then the same holds for $f(\gamma_1)$ and $f(\gamma_2)$. Indeed, if this were not true, there would be a vertical curve $\gamma$  in $DS_2$ so that the curve $f(B)$ would intersect $\gamma$ at at least two points. This is clearly impossible since vertical curves go to vertical curves under $f$ and also under its inverse. Combining this monotonicity with counting the number of slits that closed vertical curves intersect as in the proof of Theorem~\ref{T:CohopF}, we conclude that $f$ maps every closed vertical curve $\gamma$ that intersects a slit of $DS_2$ to a closed vertical curve $\gamma'$ such that $\pi(\gamma)=\pi(\gamma')$. 
This implies that $f$ preserves the distances between any two nearby slits  that intersect a closed vertical curve. 
The continuity of $f$ also gives that $f$ takes every closed vertical curve $\gamma$ to a closed vertical curve $\gamma'$ with $\pi(\gamma)=\pi(\gamma')$. 

If $\gamma$ is any vertical curve that does not intersect any slits and $p, q$ are two nearby points on $\gamma$, then there exists a sequence of closed vertical curves $\gamma_n$ that do intersect slits of $DS_2$ and such that $\gamma_n$ contains a pair of points $p_n, q_n$ with $d_{DS_2}(p,p_n), d_{DS_2}(q,q_n)\to0$ as $n\to\infty$. Moreover, the diameters of the slits intersecting $\gamma_n$ must go to zero and the points $p_n$ and $q_n$ may be chosen to be, say, the bottom-most points of the slits. Since $d_{DS_2}(p_n,q_n)=d_{DS_2}(f(p_n),f(q_n))$, taking the limit shows that $d_{DS_2}(p,q)=d_{DS_2}(f(p),f(q))$, i.e., $f$ is a local isometry on every vertical curve that does not intersect any slits. The same argument shows that if $\gamma$ is any closed vertical curve and $p,q\in\gamma$ are two nearby points  that either both do not lie on any slit, or belong to the same slit, then $d_{DS_2}(p,q)=d_{DS_2}(f(p),f(q))$. Every closed vertical curve $\gamma$ is 
locally 
a geodesic, and if it intersects slits of $DS_2$, then it can be partitioned into finitely many subcurves such that each such subcurve either does not intersect a slit  or is completely contained in one. Since $f$ restricted to every such subcurve is an isometry,  $f$ is a local isometry when restricted to every vertical curve. Thus, $f$ is an isometry on $L, R$, and its restriction to every closed vertical curve is a rotation since $f$ is an orientation preserving map on every such curve. Since $f$ fixes the endpoints  of $L$, it is the identity on it.    

From the description of its elements, it is clear that $\mathcal V$ is a normal subgroup of $\mathcal{QS}$. To show that $\mathcal{QS}=\mathcal{I}\ltimes\mathcal{V}$, it remains to show that $\mathcal{I}\cap\mathcal{V}$ is trivial. If $g$ is an isometry of $DS_2$ that is the identity on $L$,  the proof of Lemma~\ref{L:IsomDS2} shows that g is either the identity or $R_{fb}$. However, $R_{fb}$ changes the orientation of every closed vertical curve and we conclude that $\mathcal{I}\cap\mathcal{V}$ consists of the identity.
To finish the proof of the theorem we need to show that every quasisymmetry of $DS_2$ is bi-Lipschitz and $\mathcal V$ is isomorphic to $\mathcal L$. To do this, it is enough to show the latter and that every element of $\mathcal V$ is bi-Lipschitz.

Let $\tilde{DS_2}$ denote the space obtained from the infinite sequence $(S_2(k))_{k\in \Z}$ of copies of $S_2$ by gluing the top side $T(k)$ of the outer square of $S_2(k)$ to the bottom side $B(k+1)$ of the outer square of $S_2(k+1)$ using an isometry. We endow $\tilde{DS_2}$ with the path metric $d$ induced from that on $S_2$. Let $\alpha$ be a locally isometric map of $\tilde{DS_2}$ onto $DS_2$ so that $S_2(k)$ gets mapped isometrically onto the front or the back copy of $S_2$ in $DS_2$ depending on whether $k$ is even or odd. The projection $\pi$ of $S_2$ onto $\bar Q_0$ induces a 1-Lipschitz map $\tilde\pi$ of $\tilde{DS_2}$ onto the infinite strip $V=\{(x,y)\in\R^2\: 0\leq x\leq1\}$. We assume that $\tilde\pi(B(0))=\{(x,0)\in\R^2\: 0\leq x\leq1\}$.

Let $f$ be an element of $\mathcal V$.
Such a map clearly lifts to $\tilde{DS_2}$, i.e., there exists a homeomorphism $\tilde f$ of $\tilde{DS_2}$ with $f\circ\alpha=\alpha\circ\tilde f$. We may assume that $\tilde f$ is the identity on the left side $L(0)$ of $S_2(0)$. Let $\tilde i$ denote the isometry of $[0,1]$ onto the bottom side 
$B(0)$ of the outer square of $S_2(0)$ so that for the isometry $i=\alpha\circ\tilde i$ from $[0,1]$ to $B$ we have $i(0)\in L$ and $i(1)\in R$.
The properties of $f$ then readily imply that 
the map $h$ given by
$$
t\mapsto \Im(\tilde{\pi}\circ \tilde{f}\circ \tilde{i}(t)),
$$ 
where $\Im$ denotes the projection of the plane $\R^2$ onto the second coordinate, is a well-defined real-valued continuous function on $[0,1]$.

We claim that $h$ is in $\mathcal L$. 
First, assume for contradiction that $h$ is not Lipschitz. Then there exist two sequences $(x_n)$ and $(y_n)$ in $[0,1]$ such that $x_n\neq y_n$ for every $n\in\N$ and 
\begin{equation}\label{E:Fnotlip}
\frac{|h(x_n)-h(y_n)|}{|x_n-y_n|}\to\infty
\end{equation}
as $n\to\infty$. 
Since $h$ is a continuous function on the closed interval $[0,1]$, we may assume that  $(x_n)$ and $(y_n)$ are convergent sequences and $|x_n-y_n|\to0$ as $n\to\infty$. 
Let $p_n\in DS_2$ be such that it belongs to the same closed vertical curve as $i(x_n)$ and the distance between $p_n$ and $i(x_n)$ is equal to $|x_n-y_n|$. Since $f$ is quasisymmetric, there exists a homeomorphism $\eta\: [0,\infty)\to[0,\infty)$, such that
\begin{equation}\label{E:Fqs}
\frac{d_{DS_2}(f(i(x_n)),f(i(y_n)))}{d_{DS_2}(f(i(x_n)),f(p_n))}\leq\eta\bigg(\frac{d_{DS_2}(i(x_n),i(y_n))}{d_{DS_2}(i(x_n),p_n)}\bigg),\quad n\in\N.
\end{equation}
Because $i$ is an isometry, we get $d_{DS_2}(i(x_n),i(y_n))=|x_n-y_n|$, and from the choice of $p_n$ we also have $d_{DS_2}(i(x_n),p_n)=|x_n-y_n|$. Also, since $f$ is a local isometry on every closed vertical curve, 
$$
d_{DS_2}(f(i(x_n)),f(p_n))=d_{DS_2}(i(x_n),p_n)=|x_n-y_n|.
$$
The map $\alpha$ is a local isometry, and thus  
$$
d_{DS_2}(f(i(x_n)),f(i(y_n)))=d(\tilde f(\tilde{i}(x_n)),\tilde f(\tilde{i}(y_n))).
$$
Therefore, Inequality~(\ref{E:Fqs}) reduces to
\begin{equation}\label{E:Fbound}
\frac{d(\tilde f(\tilde i(x_n)),\tilde f(\tilde i(y_n)))}{|x_n-y_n|}\leq\eta(1).
\end{equation}
The map $\Im\circ\tilde\pi$ is however 1-Lipschitz, and thus
\begin{equation}\label{E:Flip}
|h(x_n)-h(y_n)|\leq d(\tilde f(\tilde i(x_n)),\tilde f(\tilde i(y_n))).
\end{equation}
Combining~(\ref{E:Fnotlip}), (\ref{E:Fbound}), and (\ref{E:Flip}), we obtain a contradiction proving that $h$ is Lipschitz. Because $f$ is the identity on $L$ and an isometry on $R$, we obtain that $h(0)=0$ and $h(1)$ is an integer. Let now $k/2^n$ be an arbitrary reduced dyadic number. By the definition of $DS_2$, a closed  vertical curve that contains $i(k/2^n)$ intersects slits of diameter $1/2^n$ and these slits are the distance $1/2^n$ apart. Since slits intersecting a closed vertical curve have to go under $f$ to the slits intersecting the same curve, we conclude that $\tilde f$ moves $\tilde i(k/2^n)$ the distance of an integer multiple of $2/2^n$, and thus  obtain that $h(k/2^n)=2l/2^n$ for some integer $l$.  This finishes the proof that $h$ is an element of $\mathcal L$.

An isomorphism from $\mathcal V$ to $\mathcal L$ is given by 
$$
f\mapsto \Im(\tilde \pi\circ \tilde f\circ \tilde i). 
$$
The fact that it is an injective homomorphism is immediate and it remains to check that it is onto. Let $h$ be an element of $\mathcal L$ and define $f\in\mathcal V$ as follows. 
Since $h(k/2^n)=2l/2^n$ for every reduced dyadic number $k/2^n$ and some integer $l$, the map $H$ of $V$ given by $H(x,y)=(x,y+h(x))$ lifts to a homeomorphism $\tilde f$ of  $\tilde{DS}_2$ such that $h(t)=\Im(\tilde\pi\circ\tilde f\circ\tilde i(t))$. 
The map $\tilde f$ clearly descends to a homeomorphism $f$ of $DS_2$, which is a rotation when restricted to every closed vertical curve. 
Since $h(0)=0$, the map $f$ is the identity on $L$, and since $h(1)$ is an integer, $f$ is an isometry on $R$. 
It remains to show that $f$ is bi-Lipschitz.  

Suppose that $h$ is $C$-Lipschitz, let $p$ and $q$ be any two points in $DS_2$, and let $\gamma$ be a shortest curve connecting them. Then $d_{DS_2}(f(p),f(q))$ is at most the length of $f(\gamma)$. 
Let $\tilde\gamma$ and 
$\tilde \gamma_f$ denote lifts in $\tilde{DS_2}$ under the map $\alpha$ of $\gamma$ and $f(\gamma)$, respectively. The lengths of $\gamma$ and $f(\gamma)$ are the Euclidean lengths of $\tilde\pi(\tilde\gamma)$ and $\tilde\pi(\tilde\gamma_f)$, respectively, since the metric $d_{DS_2}$ on $DS_2$ is the path metric induced by the Euclidean metric in the plane. 
Now, if $\tilde\pi(\tilde\gamma(t))=(x(t),y(t))$, then 
$$
\tilde\pi(\tilde\gamma_f(t))=(x(t), y(t)+h(x(t))), 
$$
and thus the Euclidean length of $\tilde\pi(\tilde\gamma_f)$ is at most $\sqrt{2}(1+C)$ times the Euclidean length of $\tilde\pi(\tilde\gamma)$.
Therefore, $d(f(p),f(q))\leq\sqrt{2}(1+C)d(p,q)$, i.e., $f$ is $\sqrt{2}(1+C)$-Lipschitz. The same argument shows that $f^{-1}$ is also $\sqrt{2}(1+C)$-Lipschitz and the theorem follows.
\qed

\end{document}